\input amstex
\documentstyle{amsppt}
\magnification=\magstep 1
\voffset-1cm
\TagsOnLeft
\topmatter
\title Selfcoincidences in higher codimensions  \endtitle
\author Ulrich Koschorke *\endauthor
\leftheadtext{Ulrich Koschorke}
\address Universit\"at Siegen,
Emmy Noether Campus, Walter-Flex-Str. 3,
D-57068 Siegen, Germany
\endaddress
\email koschorke\@mathematik.uni-siegen.de \endemail
\thanks * Supported in part within the German-Brazilian Cooperation by IB-BMBF.
\endthanks
\abstract
When can a map between manifolds be deformed away from itself? We describe a (normal bordism) obstruction which is often computable and in general much stronger than the classical primary obstruction in cohomology. In particular, it answers our question completely in a large dimension range.

As an illustration we give explicit criteria in three sample settings: projections from Stiefel manifolds to Grassmannians, sphere bundle projections and maps defined on spheres. In the first example a theorem of Becker and Schultz concerning the framed bordism class of a compact Lie group plays a central role; our approach yields also a very short geometric proof (included as an appendix) of this result.
\endabstract
\subjclass
Primary 55 M 20, 55 S 35, 57 R 90, 57 S 15
\endsubjclass
\endtopmatter

\input boxedeps.tex
\SetRokickiEPSFSpecial
\HideDisplacementBoxes

\define\id{\operatorname{id}}
\define\forg{\operatorname{forg}}
\define\const{\operatorname{const}}
\define\Hom{\operatorname{Hom}}

\document

\specialhead I.\ \ Introduction
\endspecialhead
Throughout this paper \ $M$ and $N$ \ denote smooth connected manifolds without boundary, of dimensions \ $m$ and $n$, resp., \ $M$ \ being compact. We say a map \ $f : M \longrightarrow N$ \ is {\it loose} \ (or \ $f|\wr f$ \ in the notation of Dold and Gon\c calves \cite{DG}) if \ $f$ is homotopic to some map \ $f'$ \ which has no coincidences with \ $f$, \ i.e. $f (x) \ne f' (x)$ \ for all \ $x \in M$.

\medskip
\flushpar
{\bf Problem:} \ {\it Give strong and computable criteria $($expressed in a language of algebraic topo\-logy$)$ for \ $f$ \ to be loose.}
\medskip

In this paper we present some results and examples which seem to indicate that normal bordism theory offers an appropriate language. Indeed, a careful analysis of the coincidence behaviour (of a suitable approximation) of \ $(f, f) : M \longrightarrow N \times N$ \ yields a triple \ $(C, g, \overline g)$ \ where
\roster
\item"(i)" \ $C$ \ is a smooth \ $(m - n)$--dimensional manifold (the coincidence locus);
\item"(ii)" \ $g : C \longrightarrow M$ \ is a continuous map (the inclusion); \ and
\item"(iii)" \ $\overline g$ \ is a vector bundle isomorphism which describes the stable normal bundle of \ $C$ \ in terms of the pullback \ $g^* (\varphi)$ \ of the virtual coefficient bundle \ $\varphi = f^* (TN) - TM$ \ over \ $M$.
\endroster

This leads to a well-defined looseness obstruction
$$
\omega (f) \ := \ [C, g, \overline g] \ \in \ \Omega_{m -n} (M; f^* (TN) - TM)
$$
in the normal bordism group which consists of the bordism classes of triples as above.

Quite generally, normal bordism theory deals with closed manifolds \ $C$, each equipped with a map \ $g$ \ into a given target space \ $X$ \ and -- very importantly -- with a \lq\lq twisted framing\rq\rq \ $\overline g$. \ E.g.\ trivial coefficient bundles lead to framed bordism groups \ $\Omega_*^{fr} (X)$ \ which will play a central role in the discussion below.

The Pontryagin-Thom construction makes normal bordism groups accessible to the calculation techniques of stable homotopy theory; recall that for example \ $\Omega_j (M; f^* (TN)  - TM), \ j \in \Bbb Z$, is isomorphic to the \ $(j + k + n)$-th homotopy group of the Thom space of the vector bundle \ $f^* (TN) \oplus \nu (M)$ \ over \ $M$ \ where \ $\nu (M) \oplus TM \cong M \times \Bbb R^{k+m}, \ k \gg 0$. For the computation of low-dimensional normal bordism groups via (co-)homological methods see also \cite{Ko 1}, 9.3.

\medskip\flushpar
{\bf Selfcoincidence theorem.} \ {\it Assume \ $m < 2n - 2$. Then \ $f$ \ is loose if and only if \ $\omega (f) = 0$.}
\medskip

This is our central result. In \S\ 2 below we give the proof which is based on the singularity theory for vector bundle morphisms (see \cite{Ko 1}). As a by-product we show also that if the map \ $f$ \ can be homotoped away from itself, then this can be achieved by an arbitrarily small deformation. Furthermore we obtain a formula expressing  \ $\omega (f)$ \ in terms of the Euler number \ $\chi (N)$ \ of \ $N$ \ and of the (normal bordism) degree of \ $f$. Often this makes explicit calculations possible.
\medskip

The natural Hurewicz homomorphism maps our invariant \ $\omega (f)$ \ to the Poincar\'e dual of the classical primary obstruction \ $o_n (f, f)$ \ in the (co-)homology of \ $M$ (in general with twisted coefficients; compare \cite{GJW}, theorem 3.3). This transition forgets the vector bundle isomorphism \ $\overline g$ \ nearly completely, keeping track only of the orientation information it carries. If \ $m = n$, this is no loss. However, in higher codimensions \ $m - n > 0$ \ the knowledge of \ $\overline g$ \ is usually crucial.
\medskip

\example{Example}
Consider the canonical projections
$$
p : V_{r, k} \ \longrightarrow \ G_{r, k} \ \quad \text{and} \quad \ \widetilde p : V_{r, k} \ \longrightarrow \ \widetilde G_{r, k}
$$
from the Stiefel manifold of orthonormal \ $k$--frames in \ $\Bbb R^r$ \ to the Grassmannian of (unoriented or oriented, resp.) \ $k$--planes through the origin in \ $\Bbb R^r$.

Then \ $\omega (p) = \omega (\widetilde p)$ \ lies in the framed bordism group \ $\Omega^{fr}_{d (k)} (V_{r, k})$ \ where \ $d (k) = \frac12 k (k - 1)$. In nearly all interesting cases the map \ $g$ \ from the coincidence locus \ $C$ into \ $V_{r, k}$ \ factors -- up to homotopy -- through a lower dimensional manifold so that the primary obstruction vanishes. Frequently \ $g$ \ is even nulhomotopic.
\endexample

\proclaim{Theorem} \ Assume \ $r \ge 2 k \ge 2$. Then: \ $p$ and $\widetilde p$ \ are loose if and only if
$$
0 \ = \ 2 \chi (G_{r, k})  \cdot  [SO (k)] \ \in \ \pi^S_{d (k)} \ .
$$
This condition holds e.g.\ if \ $k$ \ is even  \ or \ $k = 7 \ \text{or}\ 9$ \  or \ $\chi (G_{r, k}) \equiv 0 (12)$.
\endproclaim

Here a fascinating problem enters our discussion: to determine the order of a Lie group, when equipped with a left invariant framing and interpreted -- via the Pontryagin-Thom isomorphism -- as an element in the stable homotopy group of spheres \ $\pi^S_* \cong \Omega^{fr}_*$. Deep contributions were made e.g.\ by Atiyah and Smith \cite{AS}, Becker and Schultz \cite{BS}, Knapp \cite{Kn}, and Ossa \cite{O}, to name but a few (consult the summary of results and the references in \cite{O}). In particular, it is known that the invariantly framed special orthogonal group \ $SO (k)$ \ is nulbordant for \ $4 \le k \le 9, \ k \ne 5$ \ (cf.\ table 1 in \cite{O}) and that \ $24 [SO (k)] = 0$ \ and \ $2 [SO (2\ell)] = 0$ \ for all \ $k$ \ and \ $\ell$ \ (cf.\ \cite{O}, p.\ 315, and \cite{BS}, 4.7; for a short proof of this last claim see also our appendix).

On the other hand, the Euler number \ $\chi (G_{r, k})$ \ is easily calculated: it vanishes if \ $k \not\equiv r \equiv 0 (2)$ \ and equals \ $[r/2] \choose [k/2]$ \ otherwise (compare \cite{MS}, 6.3 and 6.4).

\proclaim{Corollary 1} \
Assume \ $r > k = 2$. Then \ $p$ and $\widetilde p$ \ are loose.
\endproclaim

\proclaim{Corollary 2} \
Assume \ $r \ge k = 3$. Then \ $p$ \ $($or, equivalently, $\widetilde p\ )$ \ is loose if and only if \ $r$ \ is even or \ $r \equiv 1 (12)$.
\endproclaim

This follows from the fact  that \ $[SO (3)] \in \pi^S_3 \cong \Bbb Z_{24}$ \ has order 12 (cf.\ \cite{AS}).

\proclaim{Corollary 3} \
Assume \ $r \ge k = 5, \ r \ne 7$. Then \ $p$ \ $($or, equivalently, \ $\widetilde p\ )$ \ is loose if and only if \ $r \not\equiv 5 (6)$.
\endproclaim

This follows since \ $[SO (5)]$ \ has order $3$ in \ $\pi^S_{10} \cong \Bbb Z_6$ \ (cf.\ \cite{O}).

The details of this example will be discussed in  \S\ 3.

\bigskip
Next consider the case when a map \ $f : M \longrightarrow N$ \ allows a section \ $s : N \longrightarrow M$ (i.e. $f \circ s = \id_N)$. Then clearly \ $f$ \ is loose if and only if \ $\id_N$ \ is -- or, equivalently, \ $\chi (N) = 0$ \ whenever \ $N$ \ is closed. In \S\ 4 we refine this simple observation in case \ $f$ is the projection of a suitable sphere bundle \ $S (\xi)$. \ Here the relative importance of  the \ $g$- and $\overline g$-data (fibre inclusion and \lq\lq twisted framing\rq\rq) in \ $\omega (f)$ \ can be studied explicitly via Gysin sequences. We obtain divisibility conditions for \ $\chi (N)$ \ in terms of the Euler class of \ $\xi$.

As a last illustration we discuss the case \ $M = S^m$ in \S\ 5. Our looseness obstruction determines (and is determined by) a group homomorphism
$$
\omega \ : \ \pi_m (N; y_0) \ \longrightarrow \ \Omega^{fr}_{m -n} \ \cong \pi^S_{m - n} \ .
$$
Thus when \ $m < 2n - 2$ \ a map \ $f : S^m \longrightarrow N$ \ is loose precisely if its homotopy class \ $[f]$ \ lies in the kernel of this homomorphism; in the case \ $N = S^n$ \ this holds if \ $2[f] = 0$ \ (when \ $n$ \ is even) and for all \ $f$ \ (when \ $n$ \ is odd) since \ $\omega = (1 + (-1)^n) \cdot E^\infty$ \ and the stable suspension \ $E^\infty$ \ is bijective here.

\medskip\flushpar
{\bf Acknowledgement.} \ It is a pleasure to thank {\it Daciberg Gon\c calves} \ for a question which raised  my interest in coincidence problems and for very stimulating discussions.

\vskip5mm
\specialhead \S\ 1. \ The coincidence invariant and the degree
\endspecialhead

Consider two maps \ $f_1, f_2 : M \longrightarrow N$.

If the resulting map \ $(f_1, f_2) : M \longrightarrow N \times N$ is smooth and transverse to the diagonal
$$
\Delta \ := \{ (y, y) \ \in \ N \times N \ | \ y \in N \}
$$
then the coincidence locus
$$
C \ := \{x \in M \ | \ f_1 (x) = f_2 (x)\} \ = \ (f_1, f_2)^{- 1} (\Delta)
\tag 1.1
$$
is a closed \ $(m - n)$--dimensional manifold canonically equipped with the following two data:
\roster
\item" " a continuous map
$$g : C \longrightarrow M \ \ \ \ \ \ \ \ \text{(namely the inclusion) \ \ ; and} \tag 1.2
$$
\item" " a stable tangent bundle isomorphism
$$\overline g \ : \ TC \oplus g^* (f^*_1 (TN)) \ \cong \ g^* (TM) \ \ \ \ \ \ \ \ \ \ \ \ \ \ \ \ \ \tag 1.3
$$
\endroster
(since the normal bundle \ $\nu (\Delta, N \times N)$ \ of \ $\Delta$ \ in \ $N \times N$ \ is canonically isomorphic to the pullback of the tangent bundle \ $TN$ \ under the first projection \ $p_1$).

If \ $f_1$ \ and \ $f_2$ \ are arbitrary continuous maps, apply the preceding construction to a smooth map \ $(f'_1, f'_2)$ \ which approximates \ $(f_1, f_2)$ \ and is transverse to \ $\Delta$. Then a (sufficiently small)  homotopy from \ $f_1$ \ to \ $f'_1$ \ determines an isomorphism \ $f_1^* (TN) \cong f^{'*}_1 (TN)$ \ which is canonical up to regular homotopy. In any case we obtain a well-defined normal bordism class
$$
\omega (f_1, f_2) \ := \ [C, g, \overline g] \ \in \ \Omega_{m -n} (M; f^*_1 (TN) - TM)
\tag 1.4
$$
which depends only on \ $f_1$ \ and on the homotopy class of \ $f_2$.

\proclaim{Proposition 1.5} \ If there exist maps \ $f'_i  :  M  \longrightarrow N$ \ which are homotopic to \ $f_i, \ i = 1, 2$, \ and such that \ $f'_1 (x) \not\equiv f'_2 (x)$ \ for all \ $x \in M$, \ then \ $\omega (f_1, f_2) = 0$.
\endproclaim

\demo{Proof}
The homotopy \ $f_1 \sim f'_1$ \ yields a nulbordism for \ $\omega (f_1, f'_2) = \omega (f_1, f_2)$.
\enddemo \hfill $\blacksquare$

Our approach also leads us to define the (normal bordism) degree of any map \ $f : M \longrightarrow N$ \ by
$$
\deg (f) \ := \ \omega (f, \text{constant map}) \ .
\tag 1.6
$$
It is represented by the inverse image \ $F$ \ of a regular value of (a smooth approximation of) \ $f$, \ together with the inclusion map and the obvious stable description of the tangent bundle \ $TF$.

\vskip5mm
\specialhead \S\ 2. \ Selfcoincidences
\endspecialhead

Given any continuous map \ $f : M \longrightarrow N$, \ we apply the previous discussion to the special case \ $f_1 = f_2 = f$. We obtain the two invariants
$$
\omega (f) := \omega (f, f), \quad \deg (f) \ \in \ \Omega_{m -n} (M; f^* (TN) - TM)
\tag 2.1
$$
(cf.\ 1.4 and 1.6), both lying in the same normal bordism group.

Any generic section \ $s$ \ of the vector bundle \ $f^* (TN)$ \ over \ $M$ \ gives rise to a map (which is homotopic to f) from \ $M$ \ to a tubular neighbourhood \ $U \cong \nu (\Delta, N \times N) \cong p^*_1 (TN)$ \ of the diagonal \ $\Delta$ \ in \ $N \times N$ \ (compare 1.3). The resulting coincidence locus, together with its normal bordism data, equals the zero set of \ $s$ (interpreted as a vector bundle homomorphism from the trivial line bundle \ $\underline{\Bbb R}$ \ to \ $f^* (TN))$, together with its singularity data (cf.\ \cite{Ko 1}). This locus consists of \ $f^{-1} \{ y_1, \dots \}$ \ if \ $s$ \ is the pullback of a generic section of \ $TN$ \ with zeroes \ $\{y_1, \dots \}$, \ which are regular values of (a smooth approximation of) \ $f$. \ In particular, if \ $N$ \ admits a nowhere zero vector field \ $v$ -- e.g.\ when \ $N$ \ is open -- then the map \ $f$ \ is loose (since it can be \lq\lq pushed slightly along $v$\rq\rq\ to get rid of all selfcoincidences). We conclude:

\proclaim{Theorem 2.2} \
Let \ $f : M^m \longrightarrow N^n$ \ be a continuous map between smooth closed connected manifolds.

Then the selfcoincidence invariant \ $\omega (f)$ \ (cf.\ $2.1$) is equal to the singularity invariant \ $\omega (\underline{\Bbb R}, f^* (TN))$ \ $($cf.\ \cite{Ko 1}, \S\ $2)$ and hence also to \ $\chi (N)\cdot \deg (f)$ \ $($cf.\ $1.6$; here \ $\chi (N)$ \ denotes the Euler number of $N)$.

Moreover, each of the following conditions implies the next one:
\roster
\item"(i)" \ $f^* (TN)$ \ allows a nowhere zero section over \ $M$;
\item"(ii)" \ $f$ can be {\rm approximated}\ by a map which has no coincidences with \ $f$ ;
\item"(iii)" \ $f$ \ is loose ;
\item"(iv)" \ there exist maps \ $f', f'' : M \longrightarrow N$ \ which have no coincidences and which are both homotopic to \ $f$ ; \ \ \ and
\item"(v)" \ $\omega (f) = 0$ .
\endroster

If \ $m < 2n - 2$, \ all these conditions are equivalent.
\endproclaim

Indeed, in this dimension range \ $\omega (\underline{\Bbb R}, f^* (TN))$ \ is the only obstruction to the existence of a monomorphism \ $\underline{\Bbb R} \ \hookrightarrow \ f^* (TN)$ \ (see theorem 3.7 in \cite{Ko 1}).

\medskip\flushpar
{\bf Special case 2.3 (codimension zero).} \ Assume \ $m = n \ge 0$. Then \ $f$ \ is loose if and only if \ $\omega (f) = 0$. Here the relevant normal bordism group \ $\Omega_0 (M; f^* (TN) - TM)$ \ is isomorphic to \ $\Bbb Z$ \ if \ $w_1 (M) = f^* (w_1 (N))$ \ and to \ $\Bbb Z_2$ \ otherwise. \ $\omega (f)$ \ counts the isolated zeroes of a generic section in \ $f^* (TN)$. We can concentrate these zeroes in a ball in \ $M$ \ and (after isotoping some of them -- if needed -- around loops where \ $w_1 (M) \ne f^* (w_1 (N))$, in order to change their signs) cancel all of them if \ $\omega (f) = 0$.

\medskip
For a very special case in higher codimensions compare \cite{DG}, 1.15.

\remark{Remark {\rm 2.4}} \
In order to understand and compute normal bordism obstructions, it is often helpful to use the natural forgetful homomorphisms
$$
\Omega_i (M; \varphi) @>{ \ \forg \ }>> \overline\Omega_i (M; \varphi) @>{ \ \mu \ }>> \ H_i (M; \widetilde{\Bbb Z}_\varphi ) \ \ \ .
$$
Here \ $\forg$ \ retains only the orientation information contained in the \ $\overline g$--components of a normal bordism class, and \ $\mu$ \ denotes the Hurewicz homomorphism to homology with coefficients which are twisted like the orientation line bundle of \ $\varphi$. \ The detailed analysis of \ $\forg$, \ given in \S\ 9 of \cite{Ko 1}, yields computing techniques which often permit to calculate obstructions in low dimensional normal bordism groups.
\endremark

\vskip5mm
\specialhead \S\ 3. \ Principal bundles
\endspecialhead

As a first example consider the projection \ $p : M \longrightarrow N$ \ of a smooth principal \ $G$--bundle (cf.\ \cite{St}, 8.1) over the closed manifold \ $N$, with \ $G$ \ a compact Lie group. A fixed choice of an orientation of \ $G$ \ at its unit element equips \ $G$ \ with a left invariant framing (which we will drop from the notation); it also yields a (stable) trivialization of the tangent bundle along the fibres of \ $p$ \ and hence of the coefficient bundle \ $p^* (TN) - TM$. \ Thus by theorem 2.2 our selfcoincidence invariant takes the form
$$
\omega (p) = \chi (N) \cdot [G, g = \text{fibre inclusion}] \ \in \ \Omega^{fr}_{m -n} (M) \ .
\tag 3.1
$$
If we concentrate on the normal bundle information -- which, in a way, represents the highest order component of this obstruction -- and neglect its $g$--part, we obtain the weaker invariant
$$
\omega' (p) \ := \ \const_* (\omega (p)) \ = \ \chi (N) [G] \ \in \ \Omega^{fr}_{m -n}
\tag 3.2
$$
which must also vanish whenever \ $p$ \ is loose. In other words, the Euler number \ $\chi (N)$ \ must be a multiple of the order of \ $[G]$ \ in \ $\Omega^{fr}_* \cong \pi^S_*$.

For \ $i \le 6$ \ the stable stem \ $\pi^S_i$ \ is generated by the class \ $[G]$ \ of some compact connected Lie group (e.g.\  $\pi^S_1 = \Bbb Z_2 \cdot [S^1]$ \ and \ $\pi^S_3 = \Bbb Z_{24} \cdot [SU (2)]$). However, this is not typical, and only the divisors of $72$ (if not of $24$) can be the order of such a class (see \cite{O}, theorem 1.1; note also Ossa's table 1).

As an illustration let us work out the details for the projections \ $p$ and $\widetilde p$ \ discussed in the example of the introduction. We may assume \ $r > k \ge 1$.

Let us first dispose of two elementary cases.

\smallskip\flushpar
{\it Case 1:} \ {\it $k = 1$. \ $p : S^{r -1} \ \longrightarrow \Bbb R P^{r -1}$ \ and \ $\widetilde p = \id_{S^{r -1}}$ \ are loose if and only if \ $r$ \ is even.}

This follows from 2.3 and 2.2.

\smallskip\flushpar
{\it Case 2:} \ {\it $k = r - 1$ \ or \ $k \equiv r -1 \not\equiv 0 (2)$: \ both \ $p$ and $\widetilde p$ \ are loose}.

Here \ $p^* (TG_{r, k}) \cong p^* (\Hom (\gamma, \gamma^\perp)) \cong \oplus^k p^* (\gamma^\perp)$ \ (cf.\ \cite{MS}, p.\ 70) has a nowhere zero section, be it for orientation reasons or since \ $G_{r, k}$ \ is odd--dimensional.

\medskip
Next recall that in the general setting \ $\dim (G_{r, k}) = k (r -k)$; the fibre dimension is given by
$$
d (k) \ := \ \dim (S)O (k) \ = \ \frac12 k (k - 1) \ \ .
\tag 3.3
$$
We have \ $p^* (TG_{r, k}) \cong \widetilde p^* (T \widetilde G_{r, k})$ \ and hence \ $\omega (p) = \omega (\widetilde p)$. According to theorem 2.2 this is the only looseness obstruction if
$$
r \ \ge \ \frac32 k \ - \ \frac 12   \ + \ \frac3k \ \ .
\tag 3.4
$$

Clearly the fibre of \ $p$ \ (or \ $\widetilde p$) over the point \ $(\Bbb R^k \subset \Bbb R^r)$ \ in the Grassmannian is \ $V_{k, k} = O (k)$ \ (or \ $SO (k)$, resp.). Also, up to homotopy the fibre inclusion \ $g$ \ factors through \ $V_{k, \ell}$ \ where \ $\ell := \max \{ 2 k - r, 0\} < k$; \ this is seen by rotating the vectors \ $v_{\ell +1}, \dots, v_k$ \ of a $k$--frame in \ $\Bbb R^k$ \ into the standard basis vectors \ $e_{k +1}, \dots, e_{2k - \ell}$ \ in \ $\Bbb R^r$. Except in situations which are already settled by the cases 1 and 2 above we see that the dimension of the intermediate manifold \ $V_{k, \ell}$ \ is strictly less (and often considerably so) than the fibre dimension \ $d (k)$ \ (cf.\ 3.3) so that the cohomological primary obstruction detects nothing.

In particular, if \ $r \ge 2k$ \ then \ $g$ \ is nulhomotopic and therefore all the information contained in the (complete!) non-selfcoincidence obstruction is already given by
$$
\omega' (\widetilde p) \ = \ 2 \cdot \chi (G_{r, k}) [SO (k)] \ \in \ \Omega^{fr}_{d (k)}
$$
(cf.\ 3.2 and 3.3). The theorem of the introduction and its corollaries follow. (If \ $(r, k) = (3, 2), (4, 3), (6, 5)$ or  $(8, 5)$ \ refer to case 2 above; if \ $(r, k) = (9, 5)$ \ the bordism class \ $a = [SO (5) \subset V_{9, 5}] \in \Omega^{fr}_{10} (V_{9, 5})$ \ lies in the image of \ $\Omega^{fr}_{10} (S^4) \cong \Bbb Z_6 \oplus \Bbb Z_2$ \ and hence \ $\omega (f) = 12 a = 0$.)

\vskip5mm
\specialhead \S\ 4. \ Sphere bundles
\endspecialhead

Let \ $\xi$ \ be a \ $(k + 1)$--dimensional real vector bundle over a closed manifold \ $N^n$. We want to study the coincidence question for the projection of the corresponding sphere bundle
$$
p  : \quad M  :=  S (\xi) \ \longrightarrow \ N \  \ .
$$

Decomposing the tangent bundle of \ $M$ into a \lq\lq horizontal\rq\rq\ and a \lq\lq vertical\rq\rq\ part, we obtain the canonical isomorphism
$$
TM \ \oplus \ \underline{\Bbb R} \ \cong \ p^* (TN) \ \oplus \ p^* (\xi) \ \ .
$$
Thus the following commuting diagram of Gysin sequences (cf.\ \cite{Sa}, 5.3 or \cite{Ko 1}, 9.20) turns out to be relevant.

$$
\EPSFxsize12.5cm
\BoxedEPSF{diag41l.eps}
$$

\vskip5mm

Here the transverse intersection homomorphism \ $\pitchfork$ \ can also be defined by applying \ $\mu \circ \forg$ \ (cf.\ 2.4) and then evaluating the (possibly twisted) Euler class \ $e (\xi)$. \ $\partial (1)$ is given by the inclusion of a typical fibre \ $S (\xi_{y_0}), \ \ y_o \in N$, \ with boundary framing induced from the compact unit ball in \ $\xi_{y_0}$; \ in other words, \ $\partial (1) = \deg (p)$ \ (cf. 1.6). Thus \ $\omega (p) = \partial (\chi (N))$ \ (cf.\ theorem 2.2) vanishes if and only if
$$
\chi (N) \ \in \ e (\xi) \left(\mu \circ \forg (\Omega_{k + 1} (N; - \xi))\right) \ \ .
\tag 4.2
$$
\medskip
We also have the successively weaker {\it necessary} \ conditions that \ $\chi (N)$ \ lies in the subgroups \ $e (\xi) (\mu (\overline\Omega_{k +1} (N; - \xi)))$ \ and \ $e (\xi) (H_{k +1} (N; \widetilde Z_\xi))$ \ (compare 2.4).
\bigskip

\example{Example 4.3} \ Let \ $\xi$ \ be an oriented real plane bundle. Then according to \cite{Ko 1}, 9.3
$$
\mu \circ \forg (\Omega_2 (N; - \xi)) \ = \ \ker (w_2 (\xi) : H_2 (N; \Bbb Z) \ \longrightarrow \ \Bbb Z_2) \ .
$$
Thus \ $\omega (p) = 0$ \ if and only if \  $\chi (N) \in e (\xi) (H_2 (N; \Bbb Z))$ \ {\it and} \ $\chi (N)$ \ is even. For all \ $n \ge 1$ \ this is also the precise condition for \ $p$ \ to be loose (if \ $n = 2$ \ it implies -- via a cohomology Gysin sequence -- that \ $e (p^* (TN)) = 0$\ ; \ therefore \ $p^* (TN)$  \ allows a nowhere vanishing section over the 2-skeleton and hence over all of \ $M$, since \ $\pi_2 (S^1) = 0)$.

As an illustration let us consider the case when \ $\xi$ \ is the \ $r$--th tensor power of the canonical complex line bundle over \ $\Bbb C P (q), \ q > 1$.  Then \ $p$ \ is loose if and only if \ $q + 1 \in r \Bbb Z = e (\xi) (H_2 (\Bbb C P (q); \Bbb Z))$ \ {\it and} \ $q$ \ is odd. This last condition is captured by normal bordism, but not by the weaker conditions (expressed in terms of oriented bordism or homology) mentioned above (cf.\ 4.2; compare also theorem 2.2 in \cite{DG}).
\endexample

\vskip5mm
\specialhead \S\ 5. \ Homotopy groups
\endspecialhead

Our last example deals with maps which are not fibre projections in general. Choose a local orientation of \ $N$ at a base point \ $y_0 \in N$. Then our looseness obstruction determines a group homomorphism
$$
\omega \ : \ \pi_m (N; y_0) \ \longrightarrow \ \Omega^{fr}_{m -n}
$$
as follows. If \ $n = 1$, \ then \ $\omega \equiv 0$. So assume \ $n \ge 2$ \ and let \ $x_0$ and $*$ denote the base point of \ $S^m$ \ and its antipode. Given \ $[f] \in \pi_m (N; y_0)$, the inclusions \ $\{ x_0\} \subset \ S^m - \{ *\} \ \subset S^m$ \ determine canonical isomorphisms (use transversality!)
$$
\Omega^{fr}_{m -n}  @>{\ \cong \ }>>  \Omega_{m-n} (S^m - \{ *\}; f^* (TN)) @>{ \ \cong \ }>> \Omega_{m -n} (S^m; f^* (TN))
$$
which we apply to the obstruction \ $\omega (f)$. Clearly, we just obtain a multiple of a similarly defined degree homomorphism (which in the case \ $N = S^n$ \ is the stable Freudenthal suspension). The relevant multiplying factor is the Euler number of \ $N$ \ (whether \ $N$ \ is closed or not).

\bigskip
\subhead\nofrills Appendix
\endsubhead

Our approach yields also a short proof of the following result which is very useful for calculations as in \S\ 3.
\medskip

\proclaim\nofrills{Theorem of Becker and Schultz} \ {\rm (cf.\ \cite{BS}, 4.5)}. \ Let \ $B$ \ be a compact connected Lie group and \ $G \subset B$ \ a proper closed subgroup. Then
$$\chi (B / G) \ \cdot \ [G] \ = \ 0 \qquad \text{in} \quad \Omega^{fr}_* \ \ .
$$
\endproclaim

\demo{Proof}
The left hand term is the (weak) selfcoincidence invariant \ $\omega' (p)$ \ of the projection \ $p : B  \ \longrightarrow \ B / G$ \ (cf. 3.2). But right multiplication with a path in \ $B$ \ from the unit to some element \ $b_0 \not\in G$, when composed with \ $p$, yields a deformation from \ $p$ \ to a map \ $p'$ \ which has no coincidences with \ $p$. Thus \ $p$ is loose and \ $\omega' (p) = 0$.

More directly: the left hand term is represented by the zero set of the pullback (under $p$) of a generic section of \ $T (B / G)$. But clearly \ $p^* (T (B / G))$ \ allows a (left invariant) section with empty zero set, and the two zero sets are framed bordant.
\enddemo

\proclaim{Corollary} \ $2 \cdot [SO (k)] = 0$ \ for all even \ $k \ge 2$.
\endproclaim

Indeed, \ $SO (k +1) / SO (k) \cong S^k$ \ has Euler number 2.

\vskip1cm

\head{References}
\endhead

\enddocument